\documentclass[10pt,a4paper,twoside,leqno]{article}
\usepackage{amsfonts,amssymb}
\usepackage{eufrak}
\usepackage{amscd}

\font\Bigtit=cmr10 scaled \magstep 4
\font\ebf=cmbx8
\font\erm=cmr8

\usepackage{graphicx}

\begin{document}

\thispagestyle{empty}

\begin{flushright}
PL ISSN 0459-6854
\end{flushright}
\vspace{0.5cm}
\centerline{\Bigtit B U L L E T I N}
\vspace{0.5cm}
\centerline{DE \ \  LA \ \  SOCI\'ET\'E \ \  DES \ \  SCIENCES \ \ ET \ \ DES \
\ \ LETTRES \ \ DE \ \ \L \'OD\'Z}
\vspace{0.3cm}
\noindent 2003\hfill Vol. LIII
\vspace{0.3cm}
\hrule
\vspace{5pt}
\noindent Recherches sur les d\'eformations \hfill Vol. XLII
\vspace{5pt}
\hrule
\vspace{0.3cm}
\noindent pp.~39--41

\vspace{0.7cm}

\noindent {\it Andrzej K. Kwa\'sniewski}

\vspace{0.7cm}

\noindent {\bf INFORMATION ON COMBINATORIAL INTERPRETATION\\ OF FIBONOMIAL
COEFFICIENTS {\small [LETTER TO THE EDITOR]}}

\vspace{0.7cm}

\noindent {\ebf Summary}

{\small   A combinatorial interpretation of the Fibonomial coefficients is
given. (Cf.: binomial Newton and Gau\ss coefficients, $q$-binomial
coefficients -- ref. [1] and refs. given therein.)}

\vspace{0.7cm}

\renewcommand{\thesubsection}{\arabic{subsection}.}

The Fibonacci sequence origin is attributed and referred to the first
edition (lost)  of  ``Liber abaci'' (1202) by Leonardo Fibonacci  [Pisano]
(see second edition from 1228 reproduced as Il Liber Abaci di Leonardo
Pisano publicato secondo la lezione Codice Maglibeciano by  Baldassarre
Boncompagni  in Scritti di Leonardo Pisano  vol. 1, (1857) Rome).

Very recently [1, 2] Fibonomial coefficients [--5] have been given
a combinatorial interpretation as counting the number of specific
finite ``birth-selfsimilar'' subposets of an infinite non-tree
poset naturally related to the Fibonacci tree of rabbits growth
process.( More: in  [2] an explicit expression for $\zeta$
characteristic function in terms of Kronecker $\delta$ of the
corresponding incidence algebra was provided- see below for
further explanatory links ).

\vspace{3mm}

In order to supply this combinatorial interpretation we  use  the
following notation: {\it  Fibonomial  coefficients} [3--5] are
defined as
$$
\left( \begin{array}{c} n\\k\end{array}
\right)_{F}=\frac{F_{n}!}{F_{k}!F_{n-k}!}\equiv
\frac{n_{F}^{\underline{k}}}{k_{F}!},\quad n_{F}\equiv F_{n}\neq 0, $$

\noindent where we make an analogy driven [6, 7]  identifications
$(n>0)$:
$$
n_{F}!\equiv n_{F}(n-1)_{F}(n-2)_{F}(n-3)_{F}\ldots 2_{F}1_{F};$$
$$0_{F}!=1;\quad n_{F}^{\underline{k}}=n_{F}(n-1)_{F}\ldots (n-k+1)_{F}. $$

\noindent This is the specification of the notation from [6] for
the purpose  Fibonomial Calculus case (see Example 2.1 in this
Bulletin paper [7]).

Define now a partially ordered infinite set  $P$ via its finite part --
subposet  $P_{m}$ ({\it rooted at $F_{1}$ level  subposet})   to be
continued ad infinitum in an obvious way as seen from the figure of
$P_{5}$ below.  It looks like the Fibonacci tree with a specific ``cobweb''.

\vspace{2mm}

\begin{center}

\includegraphics[width=75mm]{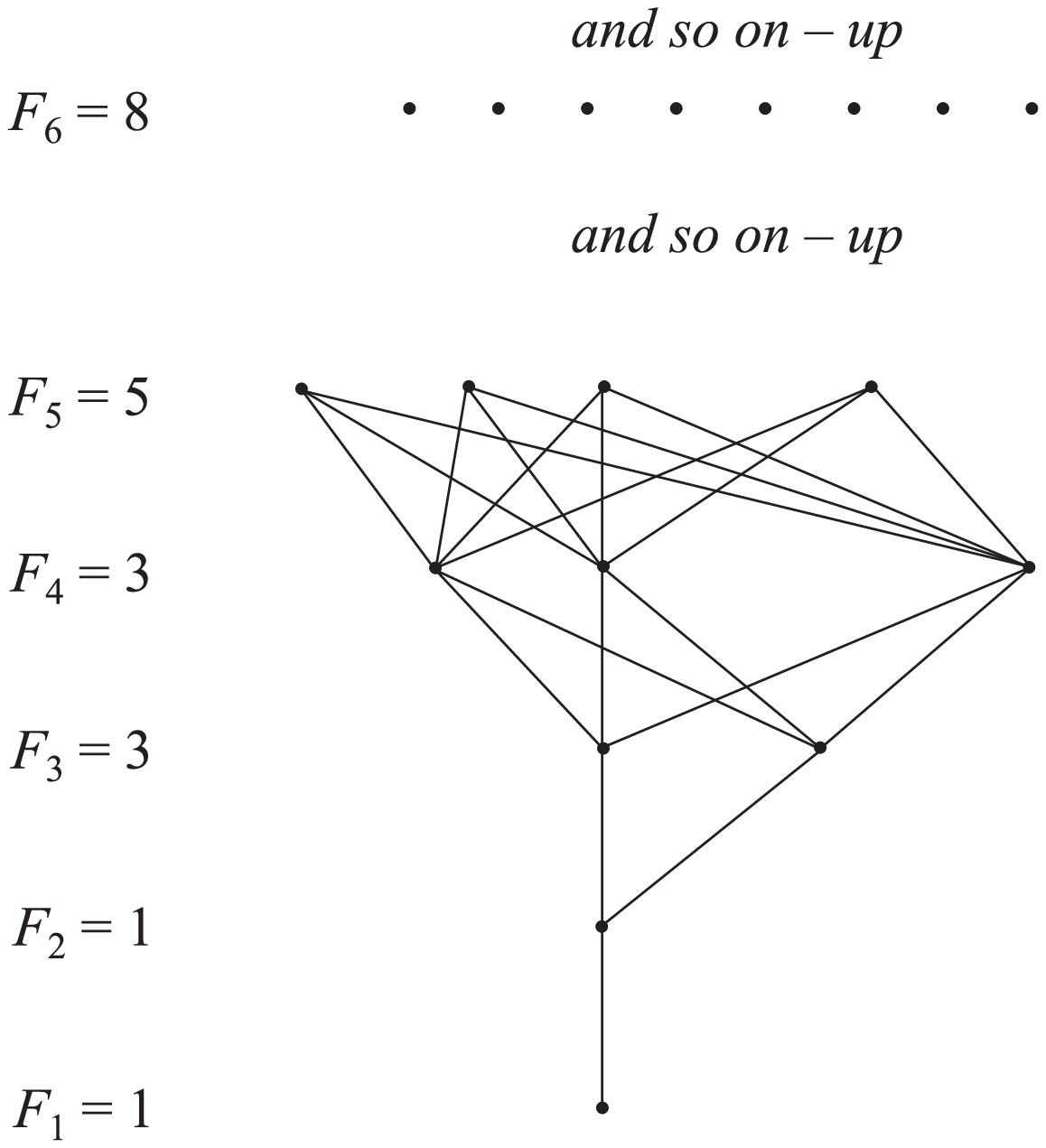}

\vspace{2mm}

\noindent {\small Fig.~1. Combinatorial interpretation of Fibonomial
coefficients.} \end{center}

\vspace{2mm}

\noindent If  one defines this poset $P$  with help of  its
incidence matrix   $\zeta$    representing P uniquely   then one
arrives at $\zeta$ with easily recognizable staircase-like
structure -- of zeros in the upper part of this upper triangle
matrix   $\zeta$; see   [1,2].

Recall:   incidence  matrix   $\zeta$ is defined for any poset as  follows
$(x,y \in P)$:
$$
  \zeta (x,y)=\left\{ \begin{array}{cl} 1&for\ x \leq y,\\0&otherwise.
\end{array} \right. $$

\vspace{2mm}

\noindent {\bf Observation 1}

{\it The number of maximal chains starting from the root  (level  $F_{1}$)
to reach any point at the n-th level  labeled by  $F_{n}$  is equal to
$n_{F}!$}.

\vspace{2mm}

\noindent {\bf Observation 2} $(k>0)$

{\it The number of maximal chains starting from any fixed  point at  the
level labeled by  $F_{k}$  to reach any point at the n-th level  labeled by
$F_{n}$ is equal to  $n_{F}^{\underline{m}}\ (n = k+m )$. }

\vspace{2mm}

\noindent {\bf Observation 3} $(k>0)$

{\it Let  $n = k+m$. The number of subposets  $P_{m}$ rooted at any fixed
point at the level labeled by  $F_{k}$  and  ending at the n-th level
labeled by $F_{n}$ is equal to}
$$ \left( \begin{array}{c} n\\k\end{array}
\right)_{F}\equiv \frac{n_{F}^{\underline{k}}}{k_{F}!}. $$

\vspace{0.5cm}

\noindent {\erm High School of Mathematics and Applied Informatics}

\noindent{\erm  Kamienna 17, PL-15-021 Bia\l ystok}

\noindent {\erm Poland}

\vspace{0.5cm}

\noindent Presented by Julian \L awrynowicz at the Session of the
Mathematical-Physical Commission of the \L \'od\'z Society of Sciences and
Arts on December 17, 2003

\vspace{0.5cm}
\noindent {\bf O INTERPRETACJI KOMBINATORYCZNEJ WSP\'O\L CZYNNIK\'OW
DWUMIENNYCH [{\small LIST DO REDAKCJI}]} \vspace{0.2cm}

\noindent {\small S t r e s z c z e n i e}

{\small Przedstawiono poszukiwan\c{a} od dziesi\c{e}cioleci
kombinatoryczn\c{a} interpretacj\c{e} wsp\'o\l czyn\-nik\'ow fibonomialnych.
(Por\'ownaj: wsp\'o\l czyniki binomialne  Newtona czy Gaussa;
wsp\'o\l \-czyn\-niki $q$-binomialne; patrz  [1] i referencje tam\.ze.)}




\begin{thebibliography}{99}
\parskip 0pt
\bibitem{}
A. K. Kwa\'sniewski, {\it Combinatorial interpretation of
Fibonomial coefficients},   Inst. Comp. Sci.  UwB/Preprint no. 52,
November 2003.
\bibitem{} A. K. Kwa\'sniewski, {\it
More on combinatorial interpretation of Fibonomial coefficients},
Inst. Comp. Sci.  UwB/Preprints no. 56, November 2003.
\bibitem{} E. Lucas, {\it Th\'eorie des fonctions num\'eriques simplement
p\'eriodiques}, American Journal of Mathematics {\bf 1} (1878),
184--240; (Translated from the French by Sidney Kravitz), Ed. D.
Lind, Fibonacci Association, 1969.

\bibitem{} G. Fonten\'e, {\it G\'en\'eralisation d`une formule connue},
Nouvelles Annales de Math\'ematiques (4) {\bf 15} (1915), 112.

\bibitem{}  H. W. Gould,   {\it The bracket function and Fonten\'e-Ward
generalized binomial coefficients with applications to Fibonomial
coefficients}, The Fibonacci Quarterly {\bf 7} (1969), 23--40.

\bibitem{}  A. K. Kwa\'sniewski,   {\it Towards  $\psi$-extension of finite
operator calculus of Rota}, Rep. Math. Phys. {\bf 47} no. 4
(2001), 305--342.  ArXiv: math.CO/0402078  2004

\bibitem{}  A. K. Kwa\'sniewski, {\it On simple characterizations of Sheffer
$\Psi$-polynomials and related propositions of the calculus of
sequences}, Bull.  Soc.  Sci.  Lettres  \L \'od\'z {\bf 52},
S\'er. Rech. D\'eform. {\bf 36} (2002), 45--65.ArXiv:
math.CO/0312397  2003.

\bibitem{} A. K. Kwa\'sniewski, {\it On duality triads}, ibid. {\bf 53}
S\'er. Rech. D\'eform. {\bf 42} (2003), 11--25.   ArXiv:
math.GM/0402260  v 1  Feb.  2004

\bibitem{} A. K. Kwa\'sniewski, {\it On Fibonomial and other triangles
versus duality triads}, ibid. {\bf 53} S\'er. Rech. D\'eform. {\bf
42} (2003), 27--37.  ArXiv:  math.GM/0402288 v 1  Feb.  2004
\end{thebibliography}
\end{document}